\documentclass[12pt]{amsart}
\usepackage{amsmath, amssymb, latexsym, amsthm}
\usepackage{mathrsfs}
\usepackage{subfigure}
\usepackage[all]{xy}

\newcommand{\ed}{\end{document}}

      \newenvironment{changemargin}[2]{\begin{list}{}{
         \setlength{\topsep}{0pt}\setlength{\leftmargin}{0pt}
         \setlength{\rightmargin}{0pt}
         \setlength{\listparindent}{\parindent}
         \setlength{\itemindent}{\parindent}
         \setlength{\parsep}{0pt plus 1pt}
         \addtolength{\leftmargin}{#1}\addtolength{\rightmargin}{#2}
         }\item }{\end{list}}

\newcommand{\fD}{\mathfrak{D}}
\newcommand{\fX}{\mathfrak{X}}
\newcommand{\Onbd}{\rmO_{\mathrm{nbd}}} %{\rmO_{\mathsf{nbd}}}
\newcommand{\Omnb}{\Omega_{\mathrm{nbd}}} %{\Omega_{\mathsf{nbd}}}
\newcommand{\od}{\mathfrak{od}}
\newcommand{\Setting}[7]{\xymatrix@R=4pt@C=7pt{#1\ar@{-}[r]&#2\ar@{-}[r]&#3\\&#4\ar@{-}[u]\\
#5\ar@{-}[uu]\ar@{-}[r] & #6\ar@{-}[u]\ar@{-}[r] & #7\ar@{-}[uu]}}

\newcommand{\Bgp}{{\Z^\N}}

\newcommand{\gp}{\mbox{-\textit{\tiny gp}}}

\newcommand{\arx}[1]{\texttt{http://arxiv.org/math/#1}}
\newcommand{\bq}{\begin{quote}}
\newcommand{\eq}{\end{quote}}
\newcommand{\cl}[1]{\overline{#1}}
\newcommand{\CH}{the Continuum Hypothesis}
\newcommand{\MA}{Martin's Axiom}

\newcommand{\Cantor}{{\{0,1\}^\N}}

\newcommand{\sr}[2]{{\txt{$#1$\\$#2$}}}
\newcommand{\N}{\mathbb{N}}
\newcommand{\NN}{{\N^{\N}}}
\newcommand{\ZN}{{\Z^{\N}}}

\newcommand{\roth}{{[\N]^{\aleph_0}}}

\newcommand{\seq}[1]{\{#1\}_{n\in\N}}
\newcommand{\setseq}[1]{\{#1 : n\in\N\}}

\newcommand{\op}{\operatorname}

\newcommand{\cI}{\mathcal{I}}
\newcommand{\cJ}{\mathcal{J}}
\newcommand{\scrA}{\mathscr{A}}
\newcommand{\scrB}{\mathscr{B}}

\newcommand{\B}{\mathrm{B}}

\newcommand{\BG}{\B_\Gamma}
\newcommand{\BL}{\B_\Lambda}
\newcommand{\BT}{\B_\Tau}

\newcommand{\BO}{\B_\Omega}

\newcommand{\CG}{C_\Gamma}

\newcommand{\CO}{C_\Omega}

\newcommand{\Tau}{\mathrm{T}}
\newcommand{\cA}{\mathcal{A}}

\newcommand{\cF}{\mathcal{F}}

\newcommand{\cM}{\mathcal{M}}
\newcommand{\cN}{\mathcal{N}}

\newcommand{\rmO}{\mathrm{O}}

\newcommand{\Q}{\mathbb{Q}}
\newcommand{\R}{\mathbb{R}}
\newcommand{\cU}{\mathcal{U}}
\newcommand{\Union}{\bigcup}
\newcommand{\cV}{\mathcal{V}}
\newcommand{\cW}{\mathcal{W}}
\newcommand{\Z}{{\mathbb Z}}

\long\def\forget#1\forgotten{}
\newcommand{\ft}{\mathfrak{t}}
\newcommand{\fb}{\mathfrak{b}}
\newcommand{\fc}{\mathfrak{c}}
\newcommand{\fd}{\mathfrak{d}}
\newcommand{\fg}{\mathfrak{g}}

\newcommand{\fu}{\mathfrak{u}}
\newcommand{\fh}{\mathfrak{h}}
\newcommand{\fp}{\mathfrak{p}}

\newcommand{\w}{\omega}
\newcommand{\x}{\times}

\newcommand{\nin}{\notin}

\newcommand{\sbst}{\subseteq}

\newcommand{\sm}{\setminus}
\newcommand{\as}{\subseteq^*}%{\let\proclaim\relax}

\newcommand{\cov}{\op{cov}}
\newcommand{\add}{\op{add}}
\newcommand{\cof}{\op{cof}}
\newcommand{\cf}{\op{cf}}
\newcommand{\non}{\op{non}}

\newtheorem{thm}{Theorem}[section]

\newtheorem{prob}[thm]{Problem}

\theoremstyle{definition}
\newtheorem{defn}[thm]{Definition}
\theoremstyle{remark}

\newtheorem{upd}{Update}
\newcommand{\update}[1]{\begin{upd}#1\end{upd}}
\newcommand{\be}{\begin{enumerate}}
\newcommand{\ee}{\end{enumerate}}
\newcommand{\bi}{\begin{itemize}}
\newcommand{\itm}{\item}
\newcommand{\ei}{\end{itemize}}
\forget
\forgotten

\newcommand{\sone}{\mathsf{S}_1}
\newcommand{\sfin}{\mathsf{S}_\mathrm{fin}}
\newcommand{\ufin}{\mathsf{U}_\mathrm{fin}}
\newcommand{\Split}{\mathsf{Split}}

\title{Selection Principles and special sets of reals}
\author{Boaz Tsaban}
\thanks{Supported by the Koshland Center for Basic Research.}
\address{Department of Mathematics,
Weizmann Institute of Science, Rehovot 76100, Israel}
\curraddr{Department of Mathematics, Bar-Ilan University, Ramat-Gan 52900, Israel}
\email{tsaban@math.biu.ac.il}
\urladdr{http://www.cs.biu.ac.il/\~{}tsaban}

\begin{document}

\begin{abstract}
We give a selection of major open problems involving selective
properties, diagonalizations, and covering properties
for sets of real numbers.

This is a revision of the version published as a chapter in the book
\textbf{Open Problems in Topology II} (E. Pearl, ed.),
Elsevier B.V., 2007, 91--108.
The present version reports solutions of some problems,
uses up-to-date notation, and has updated bibliography.

Comments and further updates would be appreciated.
\end{abstract}

\maketitle

%\tableofcontents

\section{Introduction}\label{intro}

The field of \emph{Selective Principles in Mathematics} started with
Scheepers' identification and classification of common
prototypes for selective properties appearing in classical
and modern works. For surveys of the field see \cite{LecceSurvey, KocSurv, ict}.

The main four prototypes in the field are defined as follows.
Fix a topological space $X$, and let
$\scrA$ and $\scrB$ each be a collection of covers of $X$.
Following are properties which $X$ may or may not have \cite{coc1}.
\begin{itemize}
\itm[$\binom{\scrA}{\scrB}$:] Every member of $\scrA$ has a subset which is a member of $\scrB$.
\item[$\sone(\scrA,\scrB)$:]
For each sequence $\seq{\cU_n}$ of members of $\scrA$,
there exist members $U_n\in\cU_n$, $n\in\N$, such that $\setseq{U_n}\in\scrB$.
\item[$\sfin(\scrA,\scrB)$:]
For each sequence $\seq{\cU_n}$
of members of $\scrA$, there exist finite
subsets $\cF_n\sbst\cU_n$, $n\in\N$, such that $\Union_{n\in\N}\cF_n\in\scrB$.
\item[$\ufin(\scrA,\scrB)$:]
For each sequence $\seq{\cU_n}$ of members of $\scrA$
which do not contain a finite subcover,
there exist finite subsets $\cF_n\sbst\cU_n$, $n\in\N$,
such that $\setseq{\cup\cF_n}\in\scrB$.
\end{itemize}

When $\scrA$ and $\scrB$ vary through topologically significant collections,
we obtain properties studied in various contexts by many authors.
We give some examples.

Fix a topological space $X$, and let $\rmO$ denote the collection of all
open covers of $X$.
In the case of metric spaces,
$\sfin(\rmO,\rmO)$ is the property shown by Hurewicz \cite{Hure25}
to be equivalent to Menger's basis property \cite{Menger24},
and $\sone(\rmO,\rmO)$ is Rothberger's property traditionally known
as $C''$ \cite{Roth41}, a property related to Borel's strong measure zero \cite{Borel}.

Considering special types of covers we obtain additional properties.
Henceforth, by \emph{cover of $X$} we mean a nontrivial one, i.e.,
such that $X$ itself is not a member of the cover.
An open cover $\cU$ of $X$
is an \emph{$\omega$-cover} if for each finite $F\sbst X$, there is $U\in\cU$ such that $F\subseteq U$.
$\cU$ is a \emph{$\gamma$-cover} of $X$ if it is infinite and for each $x\in X$,
$x$ is a member of all but finitely many members of $\cU$.
Let $\Omega$ and $\Gamma$ denote the collections of all open
$\omega$-covers and $\gamma$-covers of $X$, respectively.
Then $\ufin(\rmO,\Gamma)$ is the Hurewicz property \cite{Hure27},
and $\sone(\Omega,\Gamma)$ is the Gerlits-Nagy $\gamma$-property,
introduced in the context of function spaces \cite{GN}.
Additional properties of these types were studied by Arkhangel'ski\v{i},
Sakai, and others. Some of the properties are relatively new.

The field of selective principles studies the interrelations between all
properties defined by the above selective prototypes as well as similar
ones, and properties which do not fall into this category but that
can be related to properties which do.

In its broadest sense, the field (and even just its problems) cannot be
surveyed in a single book chapter.
We will restrict attention to its part dealing with
sets of real numbers.\footnote{This includes separable zero-dimensional
metric spaces, since such spaces are homeomorphic to subsets of the
irrational numbers.}
Even there, we omit several important topics.
Two of them---topological Ramsey theory and topological game theory---are
discussed in Scheepers' chapter.

While all problems we mention are about sets of real numbers, some of them
deal with sets of reals not defined by selective principles, and belong to the
more classical era of the field. Naturally, we usually mention problems we
are more familiar with.

The references we give are usually an accessible account of the problem or
related problems, but not necessarily the original source (which is usually cited
in the given reference). In fact, most of the problems have been around
much before being documented in a publication. Thus, most of the problems
posed here should be considered folklore.

The current chapter is a comprehensively revised and updated version of our earlier survey \cite{futurespm}.

\section{The Scheepers Diagram Problem}

Each of the properties mentioned in Section \ref{intro}, where
$\scrA,\scrB$ range over $\rmO,\Lambda,\Omega,\Gamma$,
is either void or equivalent to one in the following
diagram (where an arrow denotes implication)
\cite{coc1, coc2}.

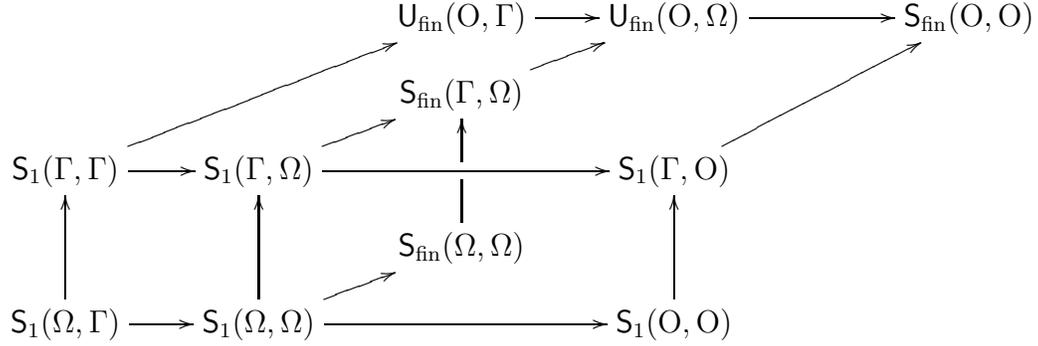
\begin{figure}[!htp]
{%\scriptsize
\renewcommand{\sr}[2]{{#1}}
\begin{changemargin}{-4cm}{-3cm}
\begin{center}
$\xymatrix@R=10pt{%@C=-2pt@R=10pt{%@=7pt{
%1
&
&
& \sr{\ufin(\rmO,\Gamma)}{\fb}\ar[r]
& \sr{\ufin(\rmO,\Omega)}{\fd}\ar[rr]
& & \sr{\sfin(\rmO,\rmO)}{\fd}
\\
%2
&
&
& \sr{\sfin(\Gamma,\Omega)}{\fd}\ar[ur]
\\
%3
& \sr{\sone(\Gamma,\Gamma)}{\fb}\ar[r]\ar[uurr]
& \sr{\sone(\Gamma,\Omega)}{\fd}\ar[rr]\ar[ur]
& & \sr{\sone(\Gamma,\rmO)}{\fd}\ar[uurr]
\\
%4
&
&
& \sr{\sfin(\Omega,\Omega)}{\fd}\ar'[u][uu]
\\
%5
& \sr{\sone(\Omega,\Gamma)}{\fp}\ar[r]\ar[uu]
& \sr{\sone(\Omega,\Omega)}{\cov(\cM)}\ar[uu]\ar[rr]\ar[ur]
& & \sr{\sone(\rmO,\rmO)}{\cov(\cM)}\ar[uu]
}$
\end{center}
\end{changemargin}
}
\caption{The Scheepers Diagram}\label{schdiag}
\end{figure}

Almost all implications which do not appear in Figure \ref{schdiag},
and are not compositions of existing implications,
are not provable: Assuming \CH{}, there are sets of reals
witnessing that \cite{coc2}.
Only the following two implications remain unsettled.

\begin{prob}[\cite{coc2}]\label{SchDiagP}
~\be
\itm Is $\ufin(\rmO,\Omega)=\sfin(\Gamma,\Omega)$?
\itm And if not, does $\ufin(\rmO,\Gamma)$ imply $\sfin(\Gamma,\Omega)$?
\ee
\end{prob}

By \emph{Borel cover} of $X$ we mean a cover of $X$ consisting of Borel subsets of $X$.
Let $\B,\BO,\BG$ denote the collections of \emph{countable Borel} covers,
$\omega$-covers, and $\gamma$-covers of $X$, respectively. Since we are
dealing with sets of reals, we may assume that all \emph{open} covers
we consider are countable \cite{split}. It follows that when $\scrA,\scrB$
range over $\B,\BO,\BG$ we get  properties stronger than the corresponding
ones when $\scrA,\scrB$ range over $\rmO,\Omega,\Gamma$. In the Borel case,
more equivalences are known and the following diagram is complete \cite{CBC}.

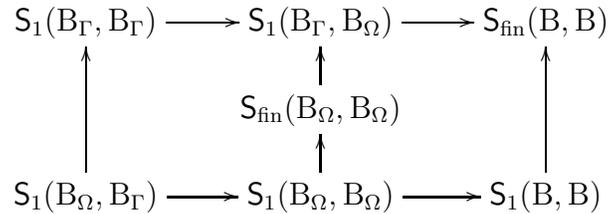
\begin{figure}[!htp]
\begin{center}
$\xymatrix@R=14pt{ %@R=10pt{%@C=-2pt@R=10pt{%@=7pt{
\sone(\BG,\BG)\ar[r] & \sone(\BG,\BO)\ar[r] & \sfin(\B,\B)\\
                     & \sfin(\BO,\BO)\ar[u]\\
\sone(\BO,\BG)\ar[r]\ar[uu] & \sone(\BO,\BO)\ar[r]\ar[u] & \sone(\B,\B)\ar[uu]
}
$
\end{center}
\caption{The Scheepers Diagram in the Borel case}\label{schdiagborel}
\end{figure}

In particular, the answer to the Borel counterpart of Problem
\ref{SchDiagP} is positive.

Problem \ref{SchDiagP} can be reformulated in terms
of topological properties of sets generating Borel non-$\sigma$-compact groups \cite{ZdBorel}.
This is related to the following problem.

\begin{prob}[\cite{ZdBorel}]\label{zdgp}
Can a Borel non-$\sigma$-compact subgroup of a Polish group be generated by a
subspace satisfying $\ufin(\rmO,\Gamma)$?
\end{prob}

%\noindent\emph{Update:} The following problem was solved. Details follow.

A set $X\sbst\R$ satisfies $\ufin(\rmO,\Gamma)$ if, and only if,
for each $G_\delta$ $G\sbst\R$ containing $X$, there is a
$\sigma$-compact $F\sbst\R$ with $X\subset F\subset G$ \cite{coc2}.

\begin{prob}[{\cite[preprint version]{BZSepar}}]\label{separ}
Assume that $X\sbst B\sbst\R$, $B$ is Borel, and $X$ satisfies $\ufin(\rmO,\Gamma)$.
Must there be a $\sigma$-compact $F$ with $X\subset F\subset B$?
What if $B$ is $F_{\sigma\delta}$?
\end{prob}

\update{``No'', and $\fb=\fd$ suffices for that \cite[final version]{BZSepar}.}

The motivation for Problem \ref{separ} was that a positive answer for its first part
would imply a negative answer for Problem \ref{zdgp}, and a positive answer for its second part would imply
a positive answer for Problem \ref{SchDiagP}(2).

\section{Examples without special set theoretic hypotheses}\label{ZFC}

\subsection{Dichotomic examples}
Let $\cJ$ be a property of sets of reals.
Sometimes there is a set theoretic hypothesis
$P$ independent of ZFC, that can be used to construct an $X\in\cJ$,
and such that its negation $\lnot P$ also implies the existence of some $Y\in\cJ$
(possibly on trivial grounds). In this case, the \emph{existence} of an
$X\in\cJ$ is a theorem of ZFC.

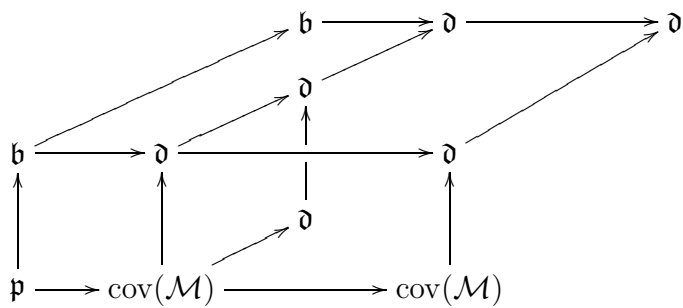
\begin{figure}[!htp]
{%\scriptsize
\renewcommand{\sr}[2]{{#2}}
\begin{changemargin}{-4cm}{-3cm}
\begin{center}
$\xymatrix@R=10pt{%@C=-2pt@R=10pt{%@=7pt{
%1
&
&
& \sr{\ufin(\rmO,\Gamma)}{\fb}\ar[r]
& \sr{\ufin(\rmO,\Omega)}{\fd}\ar[rr]
& & \sr{\sfin(\rmO,\rmO)}{\fd}
\\
%2
&
&
& \sr{\sfin(\Gamma,\Omega)}{\fd}\ar[ur]
\\
%3
& \sr{\sone(\Gamma,\Gamma)}{\fb}\ar[r]\ar[uurr]
& \sr{\sone(\Gamma,\Omega)}{\fd}\ar[rr]\ar[ur]
& & \sr{\sone(\Gamma,\rmO)}{\fd}\ar[uurr]
\\
%4
&
&
& \sr{\sfin(\Omega,\Omega)}{\fd}\ar'[u][uu]
\\
%5
& \sr{\sone(\Omega,\Gamma)}{\fp}\ar[r]\ar[uu]
& \sr{\sone(\Omega,\Omega)}{\cov(\cM)}\ar[uu]\ar[rr]\ar[ur]
& & \sr{\sone(\rmO,\rmO)}{\cov(\cM)}\ar[uu]
}$
\end{center}
\end{changemargin}
}
\caption{Critical cardinalities in the Scheepers Diagram}\label{schdiagcrit}
\end{figure}

The hypotheses used in the dichotomic arguments are
often related to combinatorial cardinal characteristics
of the continuum. See \cite{BlassHBK} for a survey of these.
The \emph{critical cardinality} of a nontrivial family $\cJ$ of sets of reals
is
$$\non(\cJ)=\min\{|X| : X\sbst\R\mbox{ and }X\nin\cJ\}.$$
Figure \ref{schdiagcrit} indicates the critical cardinalities
of the properties in the Scheepers diagram \ref{schdiag}
(we use $\cM$ for the ideal of meager sets of reals).
The critical cardinalities in the Borel case are equal to those in the open case.

Dichotomic arguments imply the existence (in ZFC) of a set
of reals $X$ satisfying $\sone(\Gamma,\Gamma)$ such that $|X|=\ft$ \cite{alpha_i},
and a set of reals satisfying  $\sfin(\Omega,\Omega)$ such that $|X|=\cf(\fd)$ \cite{SFH}.
Now, $\non(\sone(\Gamma,\Gamma))=\fb$, $\non(\sfin(\Omega,\Omega))=\fd$,
and it is consistent that $\fb>\ft$ and $\fd>\cf(\fd)$.
Thus, these existence results are not satisfactory.

\begin{prob}[\cite{ideals}]
Does there exist (in ZFC) a set of reals $X$ satisfying $\sone(\Gamma,\Gamma)$ such that $|X|=\fb$?
\end{prob}

\update{``No'', and the counter-example is Laver's model \cite{BBC}.
The answer becomes ``Yes'' if we are allowed to pick, instead of
one element from each cover, a union of two elements from each cover \cite{MHP}.}

\begin{prob}[\cite{SFH}]\label{dZFC}
Does there exist (in ZFC) a set of reals satisfying  $\sfin(\Omega,\Omega)$ such that $|X|=\fd$?
\end{prob}

\subsection{Direct constructions}
Constructions which do not appeal to a dichotomy are philosophically
much more pleasing.

There is a direct construction of a set of reals $H$ satisfying
$\ufin(\rmO,\Gamma)$ such that $|H|=\fb$ (and such that $H$ does not contain a perfect set)
\cite{BaSh01}. All finite powers of this set $H$ satisfy
$\ufin(\rmO,\Gamma)$ \cite{ideals}.
In fact, $H$ can be chosen as a subgroup or even a subfield of $\R$
\cite{o-bdd, SFH}.

\update{A stronger result with simpler proofs is given in \cite{MHP}.}

There is also a direct construction of a set of reals $M$ satisfying
$\sfin(\rmO,\rmO)$ but not $\ufin(\rmO,\Gamma)$, such that $|M|=\fd$ \cite{SFH}.

\begin{prob}[\cite{SFH}]\label{nodichoMen}
Is there a direct (non-dichotomic) construction of a set of reals $M$
satisfying $\sfin(\Omega,\Omega)$ but not $\ufin(\rmO,\Gamma)$?
\end{prob}

\subsection{The Borel case}
Let $\cJ$ be a property of sets of reals.
\emph{Borel's Conjecture for $\cJ$} is the statement ``All elements of $\cJ$ are countable''.
For all but three of the properties in the Borel case, Borel's Conjecture is consistent.

\begin{prob}[\cite{pawlikowskireclaw}]\label{BCS1BGB}
Is Borel's Conjecture for $\sfin(\B,\B)$ consistent?
\end{prob}
This is the same as asking whether it is consistent
that each uncountable set of reals can be mapped onto a dominating
subset of $\NN$ by a Borel function \cite{CBC}.
Problem \ref{BCS1BGB} is also open for $\sone(\BG,\BO)$.

\begin{prob}[Magidor]
Is Borel's Conjecture for $\sfin(\B,\B)$ equivalent to Borel's Conjecture for $\sone(\BG,\BO)$?
\end{prob}

\begin{prob}[\cite{CBC}]
Is Borel's Conjecture for $\sfin(\BO,\BO)$ consistent?
\end{prob}

\section{Examples from CH or MA}

\begin{figure}[!htp]
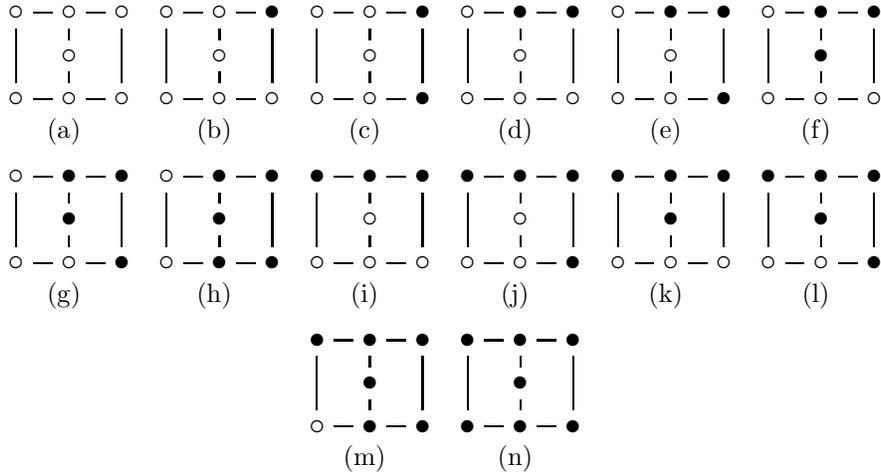

\begin{center}
\subfigure[]{$\Setting\circ\circ\circ\circ\circ\circ\circ$}
\subfigure[]{$\Setting\circ\circ\bullet\circ\circ\circ\circ$}
\subfigure[]{$\Setting\circ\circ\bullet\circ\circ\circ\bullet$}
\subfigure[]{$\Setting\circ\bullet\bullet\circ\circ\circ\circ$}
\subfigure[]{$\Setting\circ\bullet\bullet\circ\circ\circ\bullet$}
\subfigure[]{$\Setting\circ\bullet\bullet\bullet\circ\circ\circ$}
\subfigure[]{$\Setting\circ\bullet\bullet\bullet\circ\circ\bullet$}
\subfigure[]{$\Setting\circ\bullet\bullet\bullet\circ\bullet\bullet$}
\subfigure[]{$\Setting\bullet\bullet\bullet\circ\circ\circ\circ$}
\subfigure[]{$\Setting\bullet\bullet\bullet\circ\circ\circ\bullet$}
\subfigure[]{$\Setting\bullet\bullet\bullet\bullet\circ\circ\circ$}
\subfigure[]{$\Setting\bullet\bullet\bullet\bullet\circ\circ\bullet$}
\subfigure[]{$\Setting\bullet\bullet\bullet\bullet\circ\bullet\bullet$}
\subfigure[]{$\Setting\bullet\bullet\bullet\bullet\bullet\bullet\bullet$}
\end{center}
\caption{The consistent settings}\label{settings}
\end{figure}

Consider the Borel case (Figure \ref{schdiagborel}).
For each set of reals $X$, we can put ``$\bullet$'' in each place
in the diagram where the property is satisfied, and ``$\circ$'' elsewhere.
There are $14$ settings consistent with the arrows in the diagram,
and they are all listed in Figure \ref{settings}.

Setting (a) is realized by $\R\sm\Q$, i.e.\ $\NN$.

Assume \CH{}. Settings (c),(h), and (i) were realized in \cite{coc2},
Setting (k) was realized in \cite{o-bdd}, and Setting (n)
was realized in \cite{Rec94, Brendle96, MilNonGamma}.
To realize Setting (b), take a set $L$ as in Setting (c)
and a set $S$ as in Setting (i), and take $X=L\cup S$.
As $\sfin(\B,\B)$ is additive, $X$ satisfies this property.
But since $\sone(\BG,\BO)$ and $\sone(\B,\B)$ are hereditary for subsets \cite{ideals},
$X$ does not satisfy any of these.
It seems that using forcing-theoretic arguments similar to those of \cite{Brendle96},
we can realize Settings (f) and (m).

\begin{prob}
Does \CH{} imply a realization of the settings (d),(e),(g),(j), and (l)?
\end{prob}

All constructions mentioned above can be carried out using
\MA{}. Except perhaps Setting (n).

\begin{prob}[\cite{MilNonGamma}]\label{MilGa}
Does \MA{} imply the existence of a set of reals of cardinality continuum,
satisfying $\sone(\BO,\BG)$?
\end{prob}

\section{The $\delta$-property}
For a sequence $\seq{X_n}$ of subsets of $X$, define
$\liminf X_n = \Union_m\allowbreak\bigcap_{n\ge m} X_n$.
For a family $\cF$ of subsets of $X$, $L(\cF)$ denotes
its closure under the operation $\liminf$.
$X$ has the \emph{$\delta$-property} \cite{GN} if
for each open $\w$-cover $\cU$ of $X$, $X\in L(\cU)$.

Clearly, $\binom{\Omega}{\Gamma}$ implies the $\delta$-property.
$\sone(\Omega,\Gamma)=\binom{\Omega}{\Gamma}$ \cite{GN}.

\begin{prob}[\cite{GN}]\label{GNdelta}
Is the $\delta$-property equivalent to $\binom{\Omega}{\Gamma}$?
\end{prob}

\update{``Yes'' in the Borel case \cite{PCP}.}

Miller points out that, as a union of an increasing sequence of
sets with the $\delta$-property has again the $\delta$-property,
a negative answer to the following
problem implies a negative answer to Problem \ref{GNdelta}.

\begin{prob}[\cite{MilNonGamma}]\label{gammaunion}
Does every union of an increasing sequence $\seq{X_n}$ of sets satisfying $\binom{\Omega}{\Gamma}$
satisfy $\binom{\Omega}{\Gamma}$?
\end{prob}

The answer is positive in the Borel case \cite{o-bdd}.

\update{``Yes'' \cite{Jordan}.
A simplified proof and applications are given in \cite{LinSAdd}.}

\section{Preservation of properties}

\subsection{Heredity}
A property of sets of reals is \emph{hereditary} if for each set of reals $X$ satisfying
the property, all subsets of $X$ satisfy that property.
None of the selective hypotheses involving open covers is
provably hereditary \cite{ideals}. However, the property $\sone(\B,\B)$ as
well as all properties of the form $\Pi(\BG,\B)$ are hereditary
\cite{ideals} (but $\sone(\BO,\BG)$ is not \cite{MilNonGamma}).

\begin{prob}[\cite{ideals, MilNonGamma}]\label{borhered}
Is $\sone(\BO,\BO)$ or $\sfin(\BO,\BO)$ hereditary?
\end{prob}

\update{``Yes'' if each countable Borel $\omega$-cover of a set of reals
is actually an $\omega$-cover of some Borel superset of that set \cite{LinSAdd}.}

All properties in the Scheepers diagram \ref{schdiag}, except for the following
two, are known to be hereditary for $F_\sigma$ subsets \cite{SFT}.

\begin{prob}[\cite{SFT}]
Are $\sfin(\Gamma,\Omega)$ and $\sone(\Gamma,\Omega)$ hereditary for $F_\sigma$ subsets?
\end{prob}

\update{``Yes'' \cite{LinSAdd}.}

The Borel versions of all properties are hereditary for arbitrary Borel subsets \cite{CBC}.

\subsection{Finite powers}
$\sone(\Omega,\Gamma)$, $\sone(\Omega,\Omega)$, and $\sfin(\Omega,\Omega)$ are the
only properties in the open case which are preserved under taking finite powers
\cite{coc2}.

\begin{prob}[\cite{CBC}]\label{borpows}
Are any of the classes $\sone(\BO,\BO)$ or $\sfin(\BO,\BO)$ preserved by finite powers?
\end{prob}

Assume that $X$ satisfies $\sone(\BO,\BO)$ and
$Y\sbst X$. If $\sone(\BO,\BO)$ is preserved by finite powers,
then $X^k$ satisfies $\sone(\BO,\BO)$, and in particular $\sone(\B,\B)$, for all $k$.
As $\sone(\B,\B)$ is hereditary, $Y^k$ satisfies $\sone(\B,\B)$ for all $k$.
It follows that $Y$ satisfies $\sone(\BO,\BO)$ \cite{CBC}.
Similar assertions for $\sfin(\BO,\BO)$ and $\sfin(\B,\B)$ also hold \cite{CBC}.
Thus, a positive answer to Problem \ref{borpows} implies
a positive answer to Problem \ref{borhered}.

\begin{prob}[\cite{CBC}]\label{BGpow}
Is $\sone(\BO,\BG)$ preserved by finite powers?
\end{prob}

The corresponding problems for the other classes are
settled in the negative \cite{CBC}.

\subsection{Products}

Some positive results are available for products of sets.
E.g., if $X,Y\sbst\R$ have strong measure zero and
$X$ also satisfies $\ufin(\rmO,\Gamma)$,
then $X\times Y$ has strong measure zero \cite{smzpow}.

\begin{prob}[\cite{smzpow}]\label{smzp}
Assume that $X,Y\sbst\R$ satisfy $\sone(\rmO,\rmO)$, and $X$ also satisfies $\ufin(\rmO,\Gamma)$.
Does it follow that $X\times Y$ satisfies $\sone(\rmO,\rmO)$?
\end{prob}

\update{\label{prev}``No'':
\be
\itm In every extension of a model of \CH{} by $\aleph_1$ many Cohen reals,
there is a set of reals $X$ satisfying $\sone(\rmO,\rmO)$ in all finite powers and $\ufin(\rmO,\Gamma)$,
such that $X^2$ does not satisfy $\ufin(\rmO,\Gamma)$ \cite{SchTall}.
\itm The Continuum Hypothesis implies that there are sets of reals $X,Y$ satisfying $\sone(\Omega,\Gamma)$ (which is preserved by finite powers),
such that $X\times Y$ does not satisfy $\sfin(\rmO,\rmO)$ \cite{MTZ}.
\ee
}

It is not even known whether a positive answer follows when $X$ satisfies $\sone(\Omega,\Gamma)$.

\update{``No'', see Update \ref{prev}.}

The following problem withstood considerable attacks by several mathematicians.
The property in it is equivalent to the Gerlits-Nagy $(*)$ property,
and is also equivalent to $\sone(\Omega,\rmO^{\gamma\gp})$ \cite{coc7}.

\begin{prob}
Is $\ufin(\rmO,\Gamma)\cap\sone(\rmO,\rmO)$ preserved by finite products?
\end{prob}
A positive answer here implies a positive answer to Problem \ref{smzp}.
It is not even known whether $\ufin(\rmO,\Gamma)\cap\sone(\rmO,\rmO)$ is preserved by finite \emph{powers}.

\update{``No'', see Update \ref{prev}.}

None of the properties in Figure \ref{schdiag} is provably preserved by finite products
\cite{lengthdiags, CBC, huremen2, split}.
Borel's conjecture implies a consistently positive answer for $\sone(\rmO,\rmO)$ and below it.

\begin{prob}[Scheepers]\label{conpow}
Are any of the $\sfin$ or $\ufin$ type properties in the Scheepers diagram \ref{schdiag}
consistently preserved by finite products?
\end{prob}
A natural place to check Problem \ref{conpow} for $\sfin(\rmO,\rmO)$ is Miller's model
(in which, by the way, $\ufin(\rmO,\Omega)=\sfin(\rmO,\rmO)$ \cite{SF1, SFT}).

Assume that $Y$ has Hausdorff dimension zero.
The assumption that $X$ satisfies $\sone(\Omega,\Gamma)$
does not imply that $X\x Y$ has Hausdorff dimension zero.
However, if $X$ satisfies $\sone(\seq{\rmO_n},\Gamma)$,\footnote{
$\sone(\seq{\rmO_n},\Gamma)$ is the \emph{strong $\gamma$-property} \cite{GM, strongdiags}:
For each sequence $\seq{\cU_n}$ where for each $n$, $\cU_n$ is an open $n$-cover of $X$
(i.e., each $F\sbst X$ with $|F|\le n$ is contained in some member of $\cU_n$),
there are $U_n\in\cU_n$, $n\in\N$, such that $\setseq{\cU_n}$ is a $\gamma$-cover of $X$.
}
then $X\x Y$ has Hausdorff dimension zero \cite{Hdim}.

\begin{prob}[\cite{Hdim}]
Assume that $|X|<\fp$.
Is it true that for each $Y$ with Hausdorff dimension zero,
$X\x Y$ has Hausdorff dimension zero?
\end{prob}

\begin{prob}[Krawczyk]
Is it consistent (relative to ZFC) that
there are uncountable sets of reals satisfying $\sone(\Omega,\Gamma)$,
but for each such set $X$ and each set $Y$ with Hausdorff dimension zero,
$X\x Y$ has Hausdorff dimension zero?
\end{prob}

\subsection{Unions}
The question of which of the properties in Figure \ref{schdiag} is provably
preserved under taking finite or countable unions (i.e., is \emph{additive} or \emph{$\sigma$-additive})
is completely settled.
Some of the classes which are not provably additive are
\emph{consistently} additive \cite{AddQuad}.

\begin{prob}[\cite{AddQuad}]
Is $\sfin(\Omega,\Omega)$ consistently additive?
\end{prob}

The problem is also open for $\sone(\Gamma,\Omega)$ and $\sfin(\Gamma,\Omega)$.

\begin{prob}[\cite{AddQuad}]
Is $\sfin(\BO,\BO)$ consistently additive?
\end{prob}

In some cases, there remains the task to
determine the \emph{exact} additivity number.
The \emph{additivity number} of a nontrivial family $\cJ$ of sets of reals
is
$$\add(\cJ)=\min\{|\cF| : \cF\sbst\cI\mbox{ and }\cup\cF\nin\cJ\}.$$
$\max\{\fb,\fg\}\le\add(\sfin(\rmO,\rmO))\le\cf(\fd)$,
$\fh\le\add(\sone(\Gamma,\Gamma))\le\fb$,
and $\add(\cN)\le\add(\sone(\rmO,\rmO))$ \cite{AddQuad}.

\begin{prob}[\cite{AddQuad}]\label{addM}
Is $\add(\sfin(\rmO,\rmO))=\max\{\fb,\fg\}$?
\end{prob}

\begin{prob}[\cite{wqn}]\label{addGGb}
Is $\add(\sone(\Gamma,\Gamma))=\fb$?
\end{prob}

Problem \ref{addGGb} is related to Problem \ref{s1add} below.
The answer for the Borel version of Problem \ref{addGGb} is positive.

\update{Consistently, ``Yes'' \cite{BBC}.}

\begin{prob}[\cite{covM2}]
Is it consistent that $\add(\cN)<\add(\sone(\rmO,\rmO))$?
\end{prob}

Another type of problem is exemplified by the following problem.
It is easy to see that if $X$ satisfies $\sone(\Omega,\Gamma)$ and $D$ is countable,
then $X\cup D$ satisfies $\sone(\Omega,\Gamma)$.

\begin{prob}[Miller, Tsaban]
Assume that $X$ satisfies $\sone(\Omega,\Gamma)$ and $|D|<\fp$.
Does $X\cup D$ satisfy $\sone(\Omega,\Gamma)$? Is it true under \MA{} when $|D|=\aleph_1$?
\end{prob}

Recently, Francis Jordan proved that for each $D$, the following are equivalent:
\be
\itm $X\cup D$ satisfies $\sone(\Omega,\Gamma)$ for each $X$ satisfying $\sone(\Omega,\Gamma)$;
\itm $X\x D$ satisfies $\sone(\Omega,\Gamma)$ for each $X$ satisfying $\sone(\Omega,\Gamma)$.
\ee

\section{Modern types of covers}

\subsection{$\tau$-covers}\label{tau}

Recall that by ``cover of $X$'' we mean one not containing $X$ as an element.
$\cU$ is a \emph{$\large$-cover} of $X$ if each $x\in X$ is
covered by infinitely many members of $\cU$.
It is a \emph{$\tau$-cover} of $X$ if, in addition,
for each $x,y\in X$, either $\{U\in\cU : x\in U, y\nin U\}$ is finite, or
else $\{U\in\cU : y\in U, x\nin U\}$ is finite \cite{tau}.
Let $\Tau$ denote the collection of open $\tau$-covers of $X$. Then
$\Gamma \sbst \Tau \sbst \Omega$.

The most important problem concerning $\tau$-covers
is the following.

\begin{prob}[\cite{spm}]\label{omtau}
Is $\binom{\Omega}{\Gamma}=\binom{\Omega}{\Tau}$?
\end{prob}

This problem is related to many problems posed in
\cite{tau, tautau, split, strongdiags, Hdim, MilNonGamma}, etc.
The best known result in this direction is that $\binom{\Omega}{\Tau}$
implies $\sfin(\Gamma,\Tau)$ \cite{tautau}.

To state a modest form of Problem \ref{omtau}, note that if
$\binom{\Omega}{\Tau}$ implies $\sfin(\Tau,\Omega)$, then $\binom{\Omega}{\Tau}=\sfin(\Omega,\Tau)$.
\begin{prob}[\cite{tautau}]
Is $\binom{\Omega}{\Tau}=\sfin(\Omega,\Tau)$?
\end{prob}

\begin{prob}[Scheepers]\label{schtauprob}
Does $\sone(\Omega,\Tau)$ imply $\ufin(\rmO,\Gamma)$?
\end{prob}

There are many more problems of this type, and they are summarized in
\cite{MShT:858}.

Not much is known about the preservation of the new properties under set theoretic
operations.
Miller \cite{MilNonGamma} proved that assuming \CH{}, there exists
a set of reals $X$ satisfying $\sone(\BO,\BG)$ and a subset $Y$ of $X$ such that
$Y$ does not satisfy $\binom{\Omega}{\Tau}$.
Together with the remarks preceding Problem \ref{borhered}, we have that the
only classes (in addition to those in Problem \ref{borhered})
for which the heredity problem is not settled
are the following ones.

\begin{prob}[\cite{ideals}]\label{bortauhered}
Are any of the properties
$\sone(\BT,\BG)$, $\sone(\BT,\BT)$, $\sone(\BT,\BO)$, $\sone(\BT,\B)$,
$\sfin(\BT,\BT)$, or $\sfin(\BT,\BO)$, hereditary?
\end{prob}

Here are the open problems regarding unions.

\begin{prob}[\cite{AddQuad}]
Are any of the properties
$\sone(\Tau,\Tau)$, $\sfin(\Tau,\Tau)$,
$\sone(\Gamma,\Tau)$, $\sfin(\Gamma,\Tau)$,
and $\ufin(\rmO,\Tau)$ (or any of their Borel
versions) additive?
\end{prob}

It is consistent that $\ufin(\rmO,\Gamma)=\ufin(\rmO,\Tau)$,
and therefore $\ufin(\rmO,\Tau)$ is consistently $\sigma$-additive \cite{SF2}.
$\sone(\Tau,\Tau)$ is preserved under taking finite unions if, and only if,
$\sone(\Tau,\allowbreak\Tau)=\sone(\Tau,\Gamma)$ \cite{MShT:858}.

Here are the open problems regarding powers.

\begin{prob}\label{tauprod}
Are any of the properties
\be
\itm $\sone(\Omega,\Tau)$, or $\sfin(\Omega,\Tau)$,
\itm $\sone(\Tau,\Gamma)$, $\sone(\Tau,\Tau)$, $\sone(\Tau,\Omega)$,
$\sfin(\Tau,\Tau)$, or $\sfin(\Tau,\Omega)$,
\ee
preserved under taking finite powers?
\end{prob}

Most of these problems are related to Problem \ref{omtau}.

A solution to any of the problems involving $\tau$-covers
must use new ideas, since this family of covers is not as amenable as
the classical ones. In \cite{tautau} it is shown that
if we use an amenable variant of $\tau$-covers (called $\tau^*$-covers, see below),
then most of the corresponding problems can be solved.

\subsection{$\tau^*$-covers}

$Y\sbst\roth$ is \emph{linearly refinable}
if for each $y\in Y$ there exists an infinite subset
$\hat y\sbst y$ such that the family $\hat Y = \{\hat y : y\in Y\}$ is
linearly (quasi)ordered by $\as$.
A cover $\cU=\setseq{U_n}$ of $X$ is a \emph{$\tau^*$-cover} of $X$ if it is large,
and the family of all sets $\{n : x\in U_n\}$, $x\in X$,
is linearly refinable.
$\Tau^*$ is the collection of all countable open $\tau^*$-covers of $X$.

Every analytic space satisfies $\binom{\Tau}{\Gamma}$ \cite{tau}.

\begin{prob}[\cite{tautau}]
Does $\Cantor$ satisfy $\binom{\Tau^*}{\Gamma}$?
\end{prob}

\subsection{Glueable covers}\label{groupbl}

Glueability notions for covers appear naturally, under various names,
in the studies of selective principles \cite{FunRez, coc7, coc8, hureslaloms}.

A cover $\cU$ of $X$ is \emph{multifinite} if
there exists a partition of $\cU$ into infinitely many
finite covers of $X$.

\newcommand{\grp}[1]{\gimel(#1)}
\newcommand{\grpa}{\grp{\scrA}}

\begin{defn}[The \emph{Gimel operator} on families of covers]
Let $\scrA$ be a family of covers of $X$.
$\grpa$ is the family of all covers $\cU$ of $X$ such that:
Either $\cU$ is multifinite, or there exists a partition $\mathcal{P}$ of $\cU$
into finite sets such that
$\{\bigcup\mathcal{F}:\mathcal{F}\in\mathcal{P}\}\setminus\{X\}\in\scrA$.
\end{defn}

For each $\scrA$, $\scrA\sbst\grp{\scrA}$.
An element of $\grpa$ will be called \emph{$\scrA$-glueable}. This
explains our choice of the Hebrew letter \emph{Gimel} ($\gimel$).

$\grp{\Gamma}\sbst\grp{\Tau}\sbst\grp{\Omega}$.
$\sfin(\rmO,\rmO)=\sfin(\Omega,\rmO)$ \cite{coc1},
$\ufin(\rmO,\Gamma)=\sfin(\Omega,\grp{\Gamma})$ \cite{coc7},
and $\ufin(\rmO,\Omega)=\sfin(\Omega,\grp{\Omega})$ \cite{coc8}.
These results are generalized in \cite{GlCovs}.
A positive answer to the following problem is consistent \cite{SF2}.

\begin{prob}
Is $\ufin(\rmO,\Tau)=\sfin(\Omega,\grp{\Tau})$?
\end{prob}

\update{All results mentioned before the last problem are generalized in \cite{GlCovs}.
In particular, $\ufin(\rmO,\Tau^*)=\sfin(\Omega,\grp{\Tau^*})$. However, the original
problem remains open.}

$\sone(\Omega,\grp{\Omega})$ is strictly stronger than $\sone(\rmO,\rmO)$ \cite{strongdiags}.
$\sone(\Omega,\grp{\Omega})=\ufin(\rmO,\Omega)\cap\sone(\rmO,\rmO)$ \cite{strongdiags},
so the following problem can also be stated in classical terms.

\begin{prob}[\cite{coc8}]\label{s1wgp}
Is $\sone(\Omega,\Omega)=\sone(\Omega,\grp{\Omega})$?
\end{prob}

\update{$\sone(\Omega,\grp{\Omega})=\sone(\rmO,\rmO)\cap \ufin(\rmO,\Omega)$ \cite{GlCovs}.
The problem, however, is still open.}

$\ufin(\rmO,\Gamma)=\binom{\Lambda}{\grp{\Gamma}}$ \cite{hureslaloms}.
Zdomskyy proved that a positive answer to the following problem follows from NCF.

\begin{prob}
Is $\ufin(\rmO,\Omega)=\binom{\Lambda}{\grp{\Omega}}$?
\end{prob}

\section{Splittability}

Assume that $\scrA$ and $\scrB$ are collections of covers of a space $X$.
The following property was introduced in \cite{coc1}, in connection to Ramsey Theory.
\bi
\itm[$\Split(\scrA,\scrB)$:] Every cover $\cU\in\scrA$ can be split
into two disjoint subcovers $\cV$ and $\cW$, each containing an elements of $\scrB$.
\ei

If we consider this prototype with
$\scrA,\scrB\in\{\Lambda,\Omega,\Tau,\Gamma\}$, we obtain $16$
properties, each of which being either trivial or equivalent to
one in Figure \ref{splitsurv}. In this diagram, the dotted implications are open.
The implication (1) in this diagram holds if, and only if, its  implication (2) holds, and if (1) (and (2)) holds,
then (3) holds, either.

\begin{figure}[!htp]
$\xymatrix@R=16pt{
\Split(\Lambda, \Lambda) \ar[r] & \Split(\Omega, \Lambda) \ar[r] & \Split(\Tau, \Tau)\\
                           & \Split(\Omega, \Omega)\ar[u]\\
& \Split(\Omega, \Tau)\ar[u]\ar@{.>}[dr]^{(1)}\ar@/_/@{.>}[dl]_{(2)}\ar@{.>}[uul]^{(3)}\\
\Split(\Omega,\Gamma) \ar[uuu]\ar[ur]\ar[rr]     & & \Split(\Tau,\Gamma)\ar[uuu]\\
}$
\caption{}\label{splitsurv}
\end{figure}

\begin{prob}[\cite{split}]\label{splittabilityclassif}
~\be
\itm Does $\Split(\Omega, \Tau)$ imply $\Split(\Tau,\Gamma)$?
\itm And if not, then does $\Split(\Omega, \Tau)$ imply $\Split(\Lambda, \Lambda)$?
\ee
\end{prob}

The product of a $\sigma$-compact $X$ with
$Y$ satisfying $\ufin(\rmO,\scrB)$ ($\scrB\in\{\rmO,\Omega,\Gamma\}$) satisfies
$\ufin(\rmO,\scrB)$ \cite{SFT, AddQuad}.

\begin{prob}[Zdomskyy]\label{zdK}
Assume that $X$ is compact and $Y$ satisfies $\Split(\Lambda,\Lambda)$.
Does $X\x Y$ satisfy $\Split(\Lambda,\Lambda)$?
\end{prob}

Problem \ref{zdK} is also open for the other splitting properties.

\begin{prob}[\cite{AddQuad}]
Improve the lower bound or the upper bound in the inequality
$\aleph_1\le\add(\Split(\Omega, \Lambda))\le\fc$.
\end{prob}

\begin{prob}[\cite{AddQuad}]
Can the lower bound $\fu$ on $\add(\Split(\Tau, \Tau))$ be improved?
\end{prob}

All problems below are settled for the properties which do not appear in them.

\begin{prob}[\cite{split}]\label{SpLamLamAdd}
Is $\Split(\Lambda,\Lambda)$ additive?
\end{prob}

$\Split(\Lambda,\Lambda)$ is consistently additive \cite{SF1, AddQuad}.

\begin{prob}[\cite{split}]\label{BorelHered}
Are any of the properties $\Split(\BO,\BL)$, $\Split(\BO,\BO)$, $\Split(\BT,\BT)$,
and $\Split(\BT,\BG)$ hereditary?
\end{prob}

\begin{prob}[\cite{split}]\label{powcl}
Are any of the properties $\Split(\Omega,\Omega)$, $\Split(\Omega,\Tau)$, or
$\Split(\Tau,\Tau)$ preserved under taking finite powers?
\end{prob}

\section{Function spaces and local-global principles}

Let $X$ be a topological space, and $x\in X$.
A subset $A$ of $X$ \emph{converges} to $x$, $x=\lim A$,
if $A$ is infinite, $x\nin A$, and for each neighborhood $U$ of $x$, $A\sm U$ is finite.
Consider the following collections:
\begin{eqnarray*}
\Omega_x & = & \{A\sbst X : x\in\cl{A}\sm A\}\\
\Gamma_x & = & \{A\sbst X  : |A|=\aleph_0\mbox{ and }x = \lim A\}
\end{eqnarray*}
$\Gamma_x\sbst\Omega_x$. The following implications hold, and none further \cite{LG1}.
\newcommand{\arur}{\ar[u]\ar[r]}
$$\xymatrix{ %@!{
\sone(\Gamma_x,\Gamma_x)\ar[r] & \sfin(\Gamma_x,\Gamma_x)\ar[r] & \sone(\Gamma_x,\Omega_x)\ar[r] & \sfin(\Gamma_x,\Omega_x)\\
\sone(\Omega_x,\Gamma_x)\arur &\sfin(\Omega_x,\Gamma_x)\arur &\sone(\Omega_x,\Omega_x)\arur & \sfin(\Omega_x,\Omega_x)\ar[u]
}$$
In the current section, when we write $\Pi(\scrA_x,\scrB_x)$ without specifying $x$,
we mean $(\forall x)\Pi(\scrA_x,\scrB_x)$.
$\sfin(\Omega_x,\Omega_x)$ is Arkhangel'ski\v{i}'s \emph{countable fan tightness},
and $\sone(\Omega_x,\Omega_x)$ is Sakai's \emph{countable strong fan tightness}.
$\sone(\Gamma_x,\Gamma_x)$ and $\sfin(\Gamma_x,\Gamma_x)$ are
Arkhangel'ski\v{i}'s properties $\alpha_2$ and $\alpha_4$, respectively.

In the remainder of this section, $X$ will always denote a \emph{subset of $\R\sm\Q$}.
The set of all real-valued functions on $X$, denoted
$\R^X$, is equipped with the Tychonoff product topology.
$C_p(X)$ is the subspace of $\R^X$ consisting of the continuous
real-valued functions on $X$. The topology of $C_p(X)$ is known as the
\emph{topology of pointwise convergence}.
The constant zero element of $C_p(X)$ is denoted $\mathbf{0}$.

For some of the pairs $(\scrA,\scrB)\in\{\Omega,\Gamma\}^2$ and $\Pi\in\{\sone,\sfin\}$,
it is known that $C_p(X)$ satisfies $\Pi(\scrA_\mathbf{0},\scrB_\mathbf{0})$ if, and only if,
$X$ satisfies $\Pi(\scrA,\scrB)$ (see \cite{LecceSurvey} for a summary).

Fremlin's $s_1$ for $X$ and Bukovsk\'y's wQN for $X$ are
equivalent to $\sone(\Gamma_\mathbf{0},\Gamma_\mathbf{0})$ for $C_p(X)$ \cite{wqn, FrwQN}.
In a manner similar to the observation made in Section 3 of \cite{wqn},
a positive solution to Problem \ref{addGGb} should imply a positive solution
to the following problem.

\begin{prob}[\cite{FremlinSeq}]\label{s1add}
Assume that $\kappa<\fb$, and for each $\alpha<\lambda$,
$C_p(X_\alpha)$ satisfies $\sone(\Gamma_\mathbf{0},\Gamma_\mathbf{0})$.
Does $C_p(\Union_{\alpha<\kappa}X_\alpha)$ satisfy $\sone(\Gamma_\mathbf{0},\Gamma_\mathbf{0})$?
\end{prob}

\update{Consistently, ``Yes'' \cite{BBC}.}

If $X$ satisfies $\sone(\Gamma,\Gamma)$,
then $C_p(X)$ satisfies $\sone(\Gamma_\mathbf{0},\Gamma_\mathbf{0})$ \cite{wqn}.

\begin{prob}[\cite{wqn}]\label{wqn}
Is $\sone(\Gamma_\mathbf{0},\Gamma_\mathbf{0})$ for $C_p(X)$ equivalent to $\sone(\Gamma,\Gamma)$ for $X$?
\end{prob}

If the answer is positive, then Problems \ref{addGGb} and \ref{s1add} coincide.
There are several partial solutions to Problem \ref{wqn}:
First, for each $Y\sbst X$ $C_p(X)$ is $\sone(\Gamma_\mathbf{0},\Gamma_\mathbf{0})$, then $X$ satisfies
$\sone(\Gamma,\Gamma)$ \cite{Hales05}.
Second, $\sone(\Gamma_\mathbf{0},\Gamma_\mathbf{0})$ for $C_p(X)$ is equivalent to
$\cl{\mathsf{S}}_1(\Gamma,\Gamma)$ for $X$,
where $\cl{\mathsf{S}}_1$ is like $\sone$ with the following additional restriction on the
given $\gamma$-covers $\cU_n$: For each $n$, the family of closures of the elements of $\cU_{n+1}$
refines $\cU_n$ \cite{BuHa06}.
Finally, $\sone(\Gamma_\mathbf{0},\Gamma_\mathbf{0})$ for $C_p(X)$ is equivalent to $\sone(\CG,\CG)$ for $X$,
where $\CG$ is the collection of \emph{clopen} $\gamma$-covers of $X$ \cite{Sakai07}.
This reduces Problem \ref{wqn} to the question whether
$\sone(\Gamma,\Gamma)=\sone(\CG,\CG)$.

The following also seems to be open.

\begin{prob}[Scheepers]
Is $\sone(\Gamma_\mathbf{0},\Omega_\mathbf{0})$ for $C_p(X)$ equivalent to $\sone(\Gamma,\Omega)$ for $X$?
\end{prob}

$\sone(\Gamma_\mathbf{0},\Omega_\mathbf{0})$ for $C_p(X)$ is equivalent to $\sone(\CG,\CO)$ for $X$,
where $\CO$ is the collection of clopen $\omega$-covers of $X$ \cite{Sakai07}, so we really
want to know whether $\sone(\Gamma,\Omega)=\sone(\CG,\CO)$.

A topological space $Y$ is \emph{$\kappa$-Fr\'echet} if it satisfies $\binom{O(\Omega_x)}{\Gamma_x}$,
where $O(\Omega_x)$ is the family of elements of $\Omega_x$ which are open.

\begin{prob}[\cite{Sakai06}]
What is the minimal cardinality of a set $X\sbst\R$ such that $C_p(X)$ does not satisfy $\binom{O(\Omega_x)}{\Gamma_x}$?
\end{prob}

The answer is at least $\fb$ \cite{Sakai06}.

There are many additional important questions about these and related kinds of local-global principles.
Some of them are surveyed in \cite{Gru05}.

\section{Topological groups}

Let
$\Onbd$ denote the covers of $G$ of the form $\{g\cdot U : g\in G\}$,
where $U$ is a neighborhood of the unit element of $G$.
Okunev has introduced the property $\sfin(\Onbd,\rmO)$, traditionally called
\emph{$o$-boundedness} or \emph{Menger-boundedness}.
Let
$\Omnb$ denote the covers of $G$ of the form $\{F\cdot U : F\in [G]^{<\aleph_0}\}$,
where $U$ is a neighborhood of the unit element of $G$, such that for each $F\in [G]^{<\aleph_0}$,
$F\cdot U\neq G$.
Ko\v{c}inac has introduced $\sone(\Omnb,\Omega)$, $\sone(\Omnb,\Gamma)$, and $\sone(\Onbd,\Onbd)$,
traditionally called \emph{Scheepers-boundedness}, \emph{Hurewicz-boundedness}, and \emph{Rothberger-boundedness}.

The relations among these boundedness properties and their topological counterparts were
studied in many papers, see \cite{Hernandez, HRT, KMexample, o-bdd, BNS, coc11, Bab05, ZdBorel, BG},
and references therein.
In particular, the following diagram of implications is complete.
$$\xymatrix{
\sone(\Omnb,\Gamma)\ar[r] & \sone(\Omnb,\Omega)\ar[r] & \sfin(\Onbd,\rmO)\\
 & & \sone(\Onbd,\Onbd)\ar[u]
}$$

$\sfin(\Onbd,\rmO)$ is not provably preserved under cartesian products \cite{BNS, KMlinear, o-bdd}.

\begin{prob}[Tka\v{c}enko]\label{Tka}
Are there, in ZFC, groups $G,H$ satisfying $\sfin(\Onbd,\rmO)$ such that $G\x H$ does not satisfy $\sfin(\Onbd,\rmO)$?
\end{prob}

A topological group $G$ satisfies $\sone(\Omnb,\Omega)$ if, and only if, all finite powers of $G$
satisfy $\sfin(\Onbd,\rmO)$ \cite{coc11}.
Thus, the case where $G=H$ in Problem \ref{Tka} is related to the following problem.

\begin{prob}[\cite{BG}]\label{probs}
Is $\sone(\Omnb,\Omega)=\sfin(\Onbd,\rmO)$ for separable metrizable groups?
Specifically:
\be
\itm Does \CH{} imply the existence of a separable metrizable group $G$ satisfying
$\sfin(\Onbd,\rmO)$ but not $\sone(\Omnb,\allowbreak\Omega)$?
\itm Is it consistent that $\sfin(\Onbd,\rmO)=\sone(\Omnb,\Omega)$ for separable metrizable groups?
\ee
\end{prob}

\update{``Yes'' for (1): \CH{} implies that for each $k$, there is a subgroup of $\Z^\N$ such that
$G^k$ is $\sfin(\Onbd,\rmO)$ but $G^{k+1}$ is not \cite{SqMen}.}

If $G$ is \emph{analytic} and does not satisfy $\sone(\Omnb,\Gamma)$, then
$G^2$ does not satisfy $\sfin(\Onbd,\rmO)$. Thus,
for analytic groups, $\sone(\Omnb,\Gamma)=\sone(\Omnb,\Omega)$ \cite{BZSPM}.
Moreover, for analytic \emph{abelian} groups, $\sfin(\Onbd,\rmO)=\sone(\Omnb,\Gamma)$ \cite{BZSPM}.
For general analytic groups this is open.

\begin{prob}[\cite{BZSPM}]
Is there an analytic group satisfying $\sfin(\Onbd,\rmO)$ but not $\sone(\Omnb,\Gamma)$?
\end{prob}

It seems that $\ZN$ for boundedness properties of topological groups is
like $\R$ for topological and measure theoretic notions of smallness \cite{BG}.
Thus, unless otherwise indicated, all of the problems in the remainder of this
section are concerning subgroups of $\ZN$.

Say that $G\le\ZN$ is \emph{bounded} if $\{|g| : g\in G\}$ is bounded (with respect to $\le^*$),
where $|g|$ denotes the absolute value of $g$.
For subgroups of $\ZN$:
\be
\itm $G$ satisfies $\sone(\Omnb,\Gamma)$ if, and only if, $G$ is bounded \cite{Bab05}.
\itm $G$ satisfies $\sone(\Onbd,\Onbd)$ if, and only if, $G$ has strong measure zero \cite{coc11}.
\ee

\begin{prob}[\cite{BG}]
Is it consistent that there is $G\le\Bgp$
such that $G$ has strong measure zero, is unbounded, and
does not satisfy $\sfin(\rmO,\rmO)$?
\end{prob}

\begin{prob}[\cite{BG}]
Is it consistent that there is $G\le\Bgp$ such that
$G$ has strong measure zero and satisfies $\sfin(\rmO,\rmO)$,
but is unbounded and does not satisfy $\sone(\rmO,\rmO)$?
\end{prob}

\update{``Yes'' for the last two problems \cite{Weiss}. ``Yes'' under Martin's Axiom \cite{Pol}.}

Some open problems involve only the standard covering properties.
The following problem is related to Problem \ref{dZFC}.

\begin{prob}[\cite{o-bdd}]
Is there (in ZFC) a group $G\le\ZN$ of cardinality $\fd$ satisfying $\sfin(\rmO,\rmO)$?
\end{prob}

\begin{prob}[\cite{o-bdd}]\label{gammagp}
Does \CH{} imply the existence of a a group $G\le\ZN$
of cardinality $\fc$ satisfying $\sone(\BO,\BG)$, or at least $\sone(\Omega,\Gamma)$?
\end{prob}

\update{``Yes'' for $\sone(\Omega,\Gamma)$ \cite{LinSAdd}.
Possibly, a set of Sacks reals in the spirit of \cite{Brendle96}
would satisfy $\sone(\BO,\BG)$ in all finite powers and thus, by the methods of \cite{LinSAdd},
generate a group satisfying $\sone(\BO,\BG)$.}

Some approximations to Problem \ref{gammagp} are given in \cite{o-bdd}:
\CH{} implies the existence of groups satisfying $\sone(\BO,\BO)$ and of groups
satisfying $\sone(\BG,\BG)$ in all finite powers.
The answer for Problem \ref{gammagp} is positive if it is for
\ref{gammaunion}. It is also positive for the property $(\delta)$.
To get a complete positive answer, it suffices to
construct a set $X\sbst\ZN$ such that all finite powers of $X$ satisfy
$\sone(\BO,\BG)$. Thus, it suffices to have a positive answer for
Problem \ref{BGpow}.

Finally, recall Problem \ref{zdgp}, and see the other problems in
\cite{ZdBorel}.

\section{Cardinal characteristics of the continuum}

We mention here several problems in the field which are connected to
selective principles.

The main open problem in the field is the \emph{Minimal Tower Problem}.
This problem has motivated the study of $\tau$-covers.

\begin{prob}[\cite{vD}]
Is it consistent that $\fp<\ft$?
\end{prob}

Shelah is currently working on a possible positive solution to this
problem.

The study of $\tau^*$-covers, a variant of $\tau$-covers,
led to the following problem.
A family $\cF\subseteq\roth$ is \emph{linearly refinable}
if for each $A\in\cF$ there exists an infinite subset
$\hat A\subseteq A$ such that the family $\hat\cF = \{\hat A : A\in\cF\}$ is
linearly (quasi)ordered by $\as$.
$\fp^*$ is the minimal size of a centered family
in $\roth$ which is not linearly refinable.

$\fp = \min\{\fp^*,\ft\}$, and $\fp^*\le\fd$ \cite{tautau}.

\begin{prob}[\cite{tautau, ShTb768}]
Is $\fp=\fp^*$?
\end{prob}

A family $\cA\sbst (\roth)^\N$ is a \emph{$\tau$-family}
if for each $n$, $\{A(n) : A\in\cA\}$ is linearly ordered by $\as$.
A family $\cA\sbst(\roth)^\N$ is \emph{$o$-diagonalizable} if
there is $g\in\NN$ such that:
$$(\forall A\in\cA)(\exists n)\ g(n)\in A(n).$$
Let $\od$ denote the minimal cardinality of a $\tau$-family which is not $o$-diagonalizable.
$\non(\sone(\Tau,\rmO))=\od$ \cite{MShT:858}.
$\od$ is the ``tower version'' of $\cov(\cM)$: If we replace ``linearly ordered by $\as$''
by ``centered'' in the definition of $\od$, then we obtain $\cov(\cM)$.
Thus, $\cov(\cM)\le\od$.
If $\cov(\cM)=\aleph_1$, then $\od=\aleph_1$ either \cite{MShT:858}.

\begin{prob}[\cite{MShT:858}]
Is it consistent that $\cov(\cM)<\od$?
\end{prob}

Another variant of the minimal tower problem is the following.
For a cardinal number $\kappa>1$ (finite or infinite), define
$\fp_\kappa$ to be the minimal cardinality of a centered subset of $\roth$
which cannot be partitioned into less than $\kappa$ sets each having
a pseudo-intersection.

It is easy to see that $\fp=\fp_2=\fp_3=\dots=\fp_{\aleph_0}$, and $\fp_\ft=\ft$.
It turns out that $\fp = \fp_{\aleph_1}$ \cite{Pyt}.
We get a hierarchy of cardinals between $\fp$ and $\ft$:
$$\fp=\fp_{\aleph_1}\le \fp_{\aleph_2}\le ...\le \fp_\ft = \ft.$$

\begin{prob}
Is $\fp_{\fp}=\ft$?
\end{prob}

Finally, consider the following Ramsey-theoretic cardinal:
For a subset $Y$ of $\NN$ and $g\in\NN$, we say that
$g$ \emph{avoids middles} in $Y$ if:
\be
\itm for each $f\in Y$, $g\not\le^* f$;
\itm for all $f,h\in Y$ at least one of the sets
$\{n : f(n)<g(n)\le h(n)\}$ and $\{n : h(n)<g(n)\le f(n)\}$ is finite.
\ee
$\add(\fX,\fD)$ is the minimal cardinality $\kappa$ of a dominating
$Y\sbst\NN$ such that for each partition of $Y$ into $\kappa$ many
pieces, there is a piece such that no $g$ avoids middles in that piece.
This cardinal is studied in \cite{ShTb768}.

\begin{prob}[\cite{ShTb768}]
Is $\cov(\cM)\le\add(\fX,\fD)$?
\end{prob}

\section{Additional problems and other special sets of reals}

If a set of reals $X$ satisfies $\sfin(\rmO,\rmO)$,
then for each continuous image $Y$ of $X$ in $\NN$,
$Y$ is not dominating, that is,
the set $G=\{g\in\NN : (\exists f\in Y)\ g\le^* f\}$
is not equal to $\NN$ \cite{Hure25}.
In fact, $G$ satisfies $\sfin(\rmO,\rmO)$ \cite{MGD}.

If $X$ satisfies $\ufin(\rmO,\Omega)$,
then for each continuous image $Y$ of $X$ in $\NN$,
$\{g\in\NN : (\exists k)(\exists f_1,\dots,f_k\in Y)\ g\le^* \max\{f_1,\dots,f_k\}\}$
is not comeager \cite{SFT}.

\begin{prob}[\cite{SFT}]\label{Sch}
Assume that $X$ satisfies $\ufin(\rmO,\Omega)$.
Does it follow that
$$G=\{g\in\NN : (\exists k)(\exists f_1,\dots,f_k\in Y)\ g\le^* \max\{f_1,\dots,f_k\}\}$$
satisfies $\ufin(\rmO,\Omega)$?
\end{prob}

We now give a short selection of problems on special sets of reals.
See \cite{MilSpec} or the cited references for the definitions.

$X\sbst\R$ is a \emph{$\nu$-set} if for each $Y\sbst X$ which is nowhere dense in $X$,
$Y$ is countable (i.e., $X$ is Luzin relative to itself).
Every continuous image of a $\nu$-set has the property
assumed in the following problem.

\begin{prob}[\cite{BrownCox83}]
Assume that $X\sbst\R$, and for each $Y\sbst X$, $Y$ is concentrated on
a countable subset of $Y$. Does it follow that $X$ is a continuous image of
a $\nu$-set?
\end{prob}

\begin{prob}[\cite{BarJu93}]
Is it consistent that $\cov(\cM)=\aleph_1<\fc=\aleph_{\omega_1}$,
and there is a $\fc$-Luzin set (i.e., $L$ with $|L|=\fc$ and $|L\cap M|<\fc$ for
all meager $M\sbst\R$)?
\end{prob}

\begin{prob}[\cite{BarJu93}]
Assume that every strong measure zero set of reals is meager-additive.
Does Borel's Conjecture follow?
\end{prob}
The assumption in the last problem implies that $\cov(\cM)=\non\allowbreak(SMZ)<\cof(\cM)$.

If $X,Y\sbst\Cantor$ are meager-additive, then $X\x Y$ is a meager-additive
subset of $\Cantor\x\Cantor$. The same is true for null-additive subsets of $\Cantor$.
For the real line this is open.

\begin{prob}[\cite{prods}]
Assume that $X,Y\sbst\R$ are meager- (respectively, null-) additive.
Does it follow that $X\x Y$ is meager- (respectively, null-) additive?
\end{prob}

\update{``Yes'' for meager-additive \cite{Weiss2}.}

Weiss proved that every meager-additive subset of the Cantor
space, when viewed as a subset of $\R$
(where each $f\in\Cantor$ is identified with $\sum_n f(n)/2^n$),
is meager-additive (with respect to the usual addition in $\R$); and similarly for
null-additive.

\begin{prob}[\cite{prods}]
Assume that $X\sbst\R$ is meager- (respectively, null-) additive,
and $X\sbst [0,1]$.
Does it follow that $X$ is meager- (respectively, null-) additive
when viewed as a subset of $\Cantor$?
\end{prob}

\update{``Yes'' for meager-additive \cite{Weiss2}.}

For a set $H$, define $H_x=\{y: (x,y) \in H\}$.

\begin{prob}[Bartoszy\'nski]
Assume that $X\sbst\NN$ is nonmeager and $Y\sbst\NN$ is dominating.
Is there a Borel set $H\sbst\NN\x\NN$ such that
every meager set is contained in $H_x$ for some $x\in X\cup Y$?
\end{prob}

\subsubsection*{Acknowledgments}
We thank Ljubisa Ko\v{c}inac, Nadav Samet, Marion Scheepers,
and Lyubomyr Zdomskyy for their useful comments on this chapter.

%\forget

%\forgotten

\end{document}